\documentclass[12pt,a4paper]{amsart}
\usepackage{amssymb}
\usepackage{multirow}
\usepackage{graphicx}
\usepackage{amsmath,amsbsy,amsfonts,amsthm,latexsym,amsopn,amstext,amsxtra,euscript,amscd,mathrsfs}

\theoremstyle{plain}

\theoremstyle{definition}

\theoremstyle{remark}

\numberwithin{equation}{section}

\numberwithin{table}{section}

\numberwithin{figure}{section}
\input{tcilatex}

\begin{document}
\title{Positive integer powers of symmetric $(0,1)$-Heptadiagonal matrix}
\author{Murat GUBES}
\address{Department of Mathematics, Kamil Ozdag Science Faculty, Karamanoglu
Mehmetbey University, 70100 Campus, KARAMAN, TURKEY}
\email{mgubes@kmu.edu.tr}
\author{Durmus BOZKURT}
\address{Department of Mathematics, Faculty of Science, Sel\c cuk
University, 42250 Konya, Turkey}
\email{dbozkurt@selcuk.edu.tr}
\date{\today }

\begin{abstract}
In this paper, we derive a general expression for $m$th powers of symmetric $%
(0,1)$-heptadiagonal matrices with $n=3k$ order, $n\in 
\mathbb{N}
$ $(k=1,2,3,...,n/3)$.
\end{abstract}

\maketitle

\unitlength=1mm




\section{Introduction}

Let us define the $n$-by-$n$ symmetric $(0,1)$-heptadiagonal matrix $%
H=(h_{ij})$ as below: 
\begin{equation}
H=\left\{ 
\begin{array}{l}
1,\text{ for\ }\left\vert i-j\right\vert =3\text{ } \\ 
0,\text{other}%
\end{array}%
\right. \,.  \label{Hn}
\end{equation}

In literature, there are a lot of papers about matrix powers, determinants
and inverses (see [2-5], [7-8], [10] and [12-13]). Here,\ we get a general
expression of $m$th powers $(m\in 
\mathbb{N}
)$ for symmetric $(0,1)$-heptadiagonal matrix with $n=3k$ $(k=1,2,3,...,n/3)$
orders.

It is known that $m$th $(m\in 
\mathbb{N}
)$ power of a matrix $H$ is%
\begin{equation}
H^{m}=PJ^{m}P^{-1}  \label{1}
\end{equation}%
here $P$ is transforming matrix of $H$ and $J$ is jordan form of $H.$

Let we consider the following determinants

\begin{equation}
H_{n}(\alpha )=\left\vert 
\begin{array}{ccccccc}
\alpha & 0 & 0 & 1 & \cdots & 0 & 0 \\ 
0 & \alpha & 0 & 0 & 1 &  & 0 \\ 
0 & 0 & \alpha & 0 & 0 & \ddots & \vdots \\ 
1 & 0 & 0 & \alpha & \ddots & \ddots & 1 \\ 
\vdots & 1 & 0 & 0 & \ddots &  & 0 \\ 
0 &  & \ddots & \ddots & \ddots & \cdots & 0 \\ 
0 & 0 & \cdots & 1 & 0 & \cdots & \alpha%
\end{array}%
\right\vert  \label{2}
\end{equation}

and

\begin{equation}
\Delta _{n}(\alpha )=\left\vert 
\begin{array}{ccccccc}
\alpha & 1 & 0 & 0 & \cdots & 0 & 0 \\ 
1 & \alpha & 1 & 0 & \cdots &  & 0 \\ 
0 & 1 & \alpha & 1 & 0 & \ddots & \vdots \\ 
& 0 & 1 & \alpha & \ddots & \ddots & 0 \\ 
\vdots &  & 0 & 0 & \ddots & \ddots & 0 \\ 
0 &  &  & \ddots & \ddots & \cdots & 1 \\ 
0 & 0 & \cdots &  & 0 & 1 & \alpha%
\end{array}%
\right\vert \text{.}  \label{3}
\end{equation}

Using the determinant (\ref{2}), we find

\begin{equation}
\left\vert \lambda I-H\right\vert =H_{n}(\alpha )  \label{4}
\end{equation}%
where $\alpha =\lambda \in 
\mathbb{R}
$. Then, from (\ref{2}) and (\ref{3}), we write

\begin{equation}
H_{n}(\alpha )=\left( \Delta _{\frac{n}{3}}(\alpha )\right) ^{3}\text{.}
\label{5}
\end{equation}

By using definition of the $\Delta _{n}(\alpha )$ as in [2,3,4], the
recurrence relation is obtained as following

\begin{equation}
\Delta _{n}(\alpha )=\alpha \Delta _{n-1}(\alpha )-\Delta _{n-2}(\alpha )
\label{6}
\end{equation}%
where $\Delta _{0}(\alpha )=1,\Delta _{1}(\alpha )=\alpha ,\Delta
_{2}(\alpha )=\alpha ^{2}-1$.

By solving difference equation (\ref{6}) and substituting the equation into (%
\ref{5}), we get

\begin{equation}
\Delta _{\frac{n}{3}}(\alpha )=U_{\frac{n}{3}}(\frac{\alpha }{2})\text{ }
\label{7}
\end{equation}

\begin{equation}
H_{n}(\alpha )=\left( U_{\frac{n}{3}}(\frac{\alpha }{2})\right) ^{3}
\label{8}
\end{equation}%
where $U_{n}(x)$ is the $n$th degree Chebyshev polynomial of second kind
which is defined;

\begin{equation}
U_{n}(x)=\frac{\sin ((n+1)\arccos x)}{\sin (\arccos x)},-1\leq x\leq 1\text{
.}  \label{9}
\end{equation}

It's well known that all the roots of $U_{n}(x)$ are defined as follows
shown in [1], [6]

\begin{equation}
x_{nk}=\cos (\frac{k\pi }{n+1}),-1\leq x_{nk}\leq 1\text{ .}  \label{10}
\end{equation}

By considering (\ref{4}), (\ref{6})-(\ref{10}), we find eigenvalues of the
matrix $H$

\begin{equation}
\lambda _{k}=2\cos (\frac{3k\pi }{n+3}),k=\overline{1,\frac{n}{3}}\text{.}
\label{11}
\end{equation}

\section{Integer powers of $H$}

In this part of the paper, we find the matrix $P$ and $P^{-1}$ for the
expression $H=PJP^{-1}$. Secondly, we present a general expression of $H^{m}$
for $m\in 
\mathbb{N}
$ . Let we obtain the eigenvectors of matrix $H$ via linear homogeneous
system$\ $%
\begin{equation}
\left( \lambda _{j}I-H\right) x=0  \label{*}
\end{equation}

where $\lambda _{j}$ are eigenvalues of $H$ . We explicitly write down the
expression (\ref{*}) as

\begin{eqnarray}
\lambda _{j}x_{1}-x_{4} &=&0  \notag \\
\lambda _{j}x_{2}-x_{5} &=&0  \notag \\
\lambda _{j}x_{3}-x_{6} &=&0  \notag \\
&&\vdots  \label{**} \\
-x_{n-5}+\lambda _{j}x_{n-2} &=&0  \notag \\
-x_{n-4}+\lambda _{j}x_{n-1} &=&0  \notag \\
-x_{n-3}+\lambda _{j}x_{n} &=&0  \notag
\end{eqnarray}

By solving system of equations (\ref{**}), we get the eigenvectors of matrix 
$H$;

for $j=1,4,7,10,...,n-5,n-2$ and $k=1,2,...,\frac{n}{3}$

\begin{equation}
x_{jk}=U_{k-1}\left( \frac{\lambda _{\frac{j+2}{3}}}{2}\right)  \label{***}
\end{equation}

for $j=2,5,8,11,...,n-4,n-1$ and $k=1,2,...,\frac{n}{3}$

\begin{equation}
x_{jk}=U_{k-1}\left( \frac{\lambda _{\frac{j+1}{3}}}{2}\right)  \label{****}
\end{equation}

for $j=3,6,9,12,...,n-3,n$ and $k=1,2,...,\frac{n}{3}$

\begin{equation}
x_{jk}=U_{k-1}\left( \frac{\lambda _{\frac{j}{3}}}{2}\right)  \label{*****}
\end{equation}

Now, we get the expression (\ref{1}). Since, eigenvalues $\lambda _{k}$ $%
\left( k=1,2,...,\frac{n}{3}\right) $ are multiple, then each eigenvalue
corresponds triple jordan cell $\ J_{j}\left( \lambda _{k}\right) $ in the
matrix $J$. Thus, we obtain the jordan form of $H$ as

\begin{equation*}
J=diag\left( \lambda _{1},\lambda _{1},\lambda _{1},\lambda _{2},\lambda
_{2},\lambda _{2},...,\lambda _{\frac{n-3}{3}},\lambda _{\frac{n-3}{3}%
},\lambda _{\frac{n-3}{3}},\lambda _{\frac{n}{3}},\lambda _{\frac{n}{3}%
},\lambda _{\frac{n}{3}}\right)
\end{equation*}

Denoting $j$-th column of $P$ by $P_{j}(j=\left( \overline{1,n}\right) )$, $%
P=(P_{1},P_{2},P_{3},...,P_{n})$. Combining (\ref{***}),(\ref{****}) and (%
\ref{*****}), we achieve each column of $P$ as following;

\begin{equation}
P_{j}=\left( 
\begin{array}{c}
U_{0}(\frac{\lambda _{\frac{j+2}{3}}}{2}) \\ 
0 \\ 
0 \\ 
U_{1}(\frac{\lambda _{\frac{j+2}{3}}}{2}) \\ 
0 \\ 
0 \\ 
\vdots \\ 
U_{\frac{n-3}{3}}(\frac{\lambda _{\frac{j+2}{3}}}{2}) \\ 
0 \\ 
0%
\end{array}%
\right) ,j=1,7,...,n-5,P_{j}=\left( 
\begin{array}{c}
0 \\ 
U_{0}(\frac{\lambda _{\frac{j+1}{3}}}{2}) \\ 
0 \\ 
0 \\ 
U_{1}(\frac{\lambda _{\frac{j+1}{3}}}{2}) \\ 
0 \\ 
\vdots \\ 
0 \\ 
U_{\frac{n-3}{3}}(\frac{\lambda _{\frac{j+1}{3}}}{2}) \\ 
0%
\end{array}%
\right) ,j=2,8,...,n-4  \label{14}
\end{equation}

\begin{equation}
P_{j}=\left( 
\begin{array}{c}
0 \\ 
0 \\ 
U_{0}(\frac{\lambda _{\frac{j}{3}}}{2}) \\ 
0 \\ 
0 \\ 
U_{1}(\frac{\lambda _{\frac{j}{3}}}{2}) \\ 
\vdots \\ 
0 \\ 
0 \\ 
U_{\frac{n-3}{3}}(\frac{\lambda _{\frac{j}{3}}}{2})%
\end{array}%
\right) ,j=3,9,...,n-3,\ P_{j}=\left( 
\begin{array}{c}
U_{\frac{n-3}{3}}(\frac{\lambda _{\frac{j+2}{3}}}{2}) \\ 
0 \\ 
0 \\ 
U_{\frac{n-6}{3}}(\frac{\lambda _{\frac{j+2}{3}}}{2}) \\ 
0 \\ 
0 \\ 
\vdots \\ 
U_{0}(\frac{\lambda _{\frac{j+2}{3}}}{2}) \\ 
0 \\ 
0%
\end{array}%
\right) ,j=4,10,...,n-2  \label{15}
\end{equation}

\begin{equation*}
P_{j}=\left( 
\begin{array}{c}
0 \\ 
U_{\frac{n-3}{3}}(\frac{\lambda _{\frac{j+1}{3}}}{2}) \\ 
0 \\ 
0 \\ 
U_{\frac{n-6}{3}}(\frac{\lambda _{\frac{j+1}{3}}}{2}) \\ 
0 \\ 
\vdots \\ 
0 \\ 
U_{0}(\frac{\lambda _{\frac{j+1}{3}}}{2}) \\ 
0%
\end{array}%
\right) ,j=5,11,...,n-1,P_{j}=\left( 
\begin{array}{c}
0 \\ 
0 \\ 
U_{\frac{n-3}{3}}(\frac{\lambda _{\frac{j}{3}}}{2}) \\ 
0 \\ 
0 \\ 
U_{\frac{n-6}{3}}(\frac{\lambda _{\frac{j}{3}}}{2}) \\ 
\vdots \\ 
0 \\ 
0 \\ 
U_{0}(\frac{\lambda _{\frac{j}{3}}}{2})%
\end{array}%
\right) ,j=6,12,...n
\end{equation*}%
where $U_{k}(\lambda _{k})$ denote the second kind Chebyshev polynomials.
Hence, the transforming matrix is found as below; \ for $%
n=3k,(k=1,3,5,...,n/3)$;

\begin{eqnarray}
{\small P} &{\small =}&\left[ 
\begin{array}{ccccccc}
U_{0}(\frac{\lambda _{1}}{2}) & 0 & 0 & U_{\frac{n-3}{3}}(\frac{\lambda _{2}%
}{2}) & 0 & 0 & \cdots \\ 
0 & U_{0}(\frac{\lambda _{1}}{2}) & 0 & 0 & U_{\frac{n-3}{3}}(\frac{\lambda
_{2}}{2}) & 0 & \cdots \\ 
0 & 0 & U_{0}(\frac{\lambda _{1}}{2}) & 0 & 0 & U_{\frac{n-3}{3}}(\frac{%
\lambda _{2}}{2}) & \cdots \\ 
U_{1}(\frac{\lambda _{1}}{2}) & 0 & 0 & U_{\frac{n-6}{3}}(\frac{\lambda _{2}%
}{2}) & 0 & 0 & \ddots \\ 
0 & U_{1}(\frac{\lambda _{1}}{2}) & 0 & 0 & U_{\frac{n-6}{3}}(\frac{\lambda
_{2}}{2}) & 0 & \cdots \\ 
0 & 0 & U_{1}(\frac{\lambda _{1}}{2}) & 0 & 0 & U_{\frac{n-6}{3}}(\frac{%
\lambda _{2}}{2}) & \cdots \\ 
\vdots & \vdots & \vdots & \vdots & \vdots & \vdots & \ddots \\ 
U_{\frac{n-3}{3}}(\frac{\lambda _{1}}{2}) & 0 & 0 & U_{0}(\frac{\lambda _{2}%
}{2}) & 0 & 0 & \cdots \\ 
0 & U_{\frac{n-3}{3}}(\frac{\lambda _{1}}{2}) & 0 & 0 & U_{0}(\frac{\lambda
_{2}}{2}) & 0 & \cdots \\ 
0 & 0 & U_{\frac{n-3}{3}}(\frac{\lambda _{1}}{2}) & 0 & 0 & U_{0}(\frac{%
\lambda _{2}}{2}) & \cdots%
\end{array}%
\right.  \notag \\
&&\left. 
\begin{array}{cccc}
\cdots & U_{0}(\frac{\lambda _{\frac{n}{3}}}{2}) & 0 & 0 \\ 
\cdots & 0 & U_{0}(\frac{\lambda _{\frac{n}{3}}}{2}) & 0 \\ 
\cdots & 0 & 0 & U_{0}(\frac{\lambda _{\frac{n}{3}}}{2}) \\ 
\ddots & U_{1}(\frac{\lambda _{\frac{n}{3}}}{2}) & 0 & 0 \\ 
\cdots & 0 & U_{1}(\frac{\lambda _{\frac{n}{3}}}{2}) & 0 \\ 
\cdots & 0 & 0 & U_{1}(\frac{\lambda _{\frac{n}{3}}}{2}) \\ 
\ddots & \vdots & \vdots & \vdots \\ 
\cdots & U_{\frac{n-3}{3}}(\frac{\lambda _{\frac{n}{3}}}{2}) & 0 & 0 \\ 
\cdots & 0 & U_{\frac{n-3}{3}}(\frac{\lambda _{\frac{n}{3}}}{2}) & 0 \\ 
\cdots & 0 & 0 & U_{\frac{n-3}{3}}(\frac{\lambda _{\frac{n}{3}}}{2})%
\end{array}%
\right]  \label{+}
\end{eqnarray}%
for $n=3k,(k=2,4,6,...,n/3)$;

\begin{eqnarray}
{\normalsize P} &{\normalsize =}&\left[ 
\begin{array}{ccccccc}
U_{0}(\frac{\lambda _{1}}{2}) & 0 & 0 & U_{\frac{n-3}{3}}(\frac{\lambda _{2}%
}{2}) & 0 & 0 & \cdots \\ 
0 & U_{0}(\frac{\lambda _{1}}{2}) & 0 & 0 & U_{\frac{n-3}{3}}(\frac{\lambda
_{2}}{2}) & 0 & \cdots \\ 
0 & 0 & U_{0}(\frac{\lambda _{1}}{2}) & 0 & 0 & U_{\frac{n-3}{3}}(\frac{%
\lambda _{2}}{2}) & \cdots \\ 
U_{1}(\frac{\lambda _{1}}{2}) & 0 & 0 & U_{\frac{n-6}{3}}(\frac{\lambda _{2}%
}{2}) & 0 & 0 & \ddots \\ 
0 & U_{1}(\frac{\lambda _{1}}{2}) & 0 & 0 & U_{\frac{n-6}{3}}(\frac{\lambda
_{2}}{2}) & 0 & \cdots \\ 
0 & 0 & U_{1}(\frac{\lambda _{1}}{2}) & 0 & 0 & U_{\frac{n-6}{3}}(\frac{%
\lambda _{2}}{2}) & \cdots \\ 
\vdots & \vdots & \vdots & \vdots & \vdots & \vdots & \ddots \\ 
U_{\frac{n-3}{3}}(\frac{\lambda _{1}}{2}) & 0 & 0 & U_{0}(\frac{\lambda _{2}%
}{2}) & 0 & 0 & \cdots \\ 
0 & U_{\frac{n-3}{3}}(\frac{\lambda _{1}}{2}) & 0 & 0 & U_{0}(\frac{\lambda
_{2}}{2}) & 0 & \cdots \\ 
0 & 0 & U_{\frac{n-3}{3}}(\frac{\lambda _{1}}{2}) & 0 & 0 & U_{0}(\frac{%
\lambda _{2}}{2}) & \cdots%
\end{array}%
\right.  \notag \\
&&\left. 
\begin{array}{cccc}
\cdots & U_{\frac{n-3}{3}}(\frac{\lambda _{\frac{n}{3}}}{2}) & 0 & 0 \\ 
\cdots & 0 & U_{\frac{n-3}{3}}(\frac{\lambda _{\frac{n}{3}}}{2}) & 0 \\ 
\cdots & 0 & 0 & U_{\frac{n-3}{3}}(\frac{\lambda _{\frac{n}{3}}}{2}) \\ 
\ddots & U_{\frac{n-6}{3}}(\frac{\lambda _{\frac{n}{3}}}{2}) & 0 & 0 \\ 
\cdots & 0 & U_{\frac{n-6}{3}}(\frac{\lambda _{\frac{n}{3}}}{2}) & 0 \\ 
\cdots & 0 & 0 & U_{\frac{n-6}{3}}(\frac{\lambda _{\frac{n}{3}}}{2}) \\ 
\ddots & \vdots & \vdots & \vdots \\ 
\cdots & U_{0}(\frac{\lambda _{\frac{n}{3}}}{2}) & 0 & 0 \\ 
\cdots & 0 & U_{0}(\frac{\lambda _{\frac{n}{3}}}{2}) & 0 \\ 
\cdots & 0 & 0 & U_{0}(\frac{\lambda _{\frac{n}{3}}}{2})%
\end{array}%
\right] \text{.}  \label{++}
\end{eqnarray}

Denoting $j$-th column of $P^{-1}$ by $\rho _{j}(j=\left( \overline{1,n}%
\right) )$. We obtain the $P^{-1}=(\rho _{1},\rho _{2},\rho _{3},...,\rho
_{n-2},\rho _{n-1},\rho _{n})$ matrix such that;

\begin{equation}
\rho _{j}=\left( 
\begin{array}{c}
h_{1}U_{\frac{j-1}{3}}(\frac{\lambda _{1}}{2}) \\ 
0 \\ 
0 \\ 
h_{2}U_{\frac{j-1}{3}}(\frac{\lambda _{2}}{2}) \\ 
0 \\ 
0 \\ 
\vdots \\ 
h_{\frac{n}{3}}U_{\frac{j-1}{3}}(\frac{\lambda _{\frac{n}{3}}}{2}) \\ 
0 \\ 
0%
\end{array}%
\right) ,j=1,4,7,...,n-2\text{ , }\rho _{j}=\left( 
\begin{array}{c}
0 \\ 
h_{1}U_{\frac{j-2}{3}}(\frac{\lambda _{1}}{2}) \\ 
0 \\ 
0 \\ 
h_{2}U_{\frac{j-2}{3}}(\frac{\lambda _{2}}{2}) \\ 
0 \\ 
\vdots \\ 
0 \\ 
h_{\frac{n}{3}}U_{\frac{j-2}{3}}(\frac{\lambda _{\frac{n}{3}}}{2}) \\ 
0%
\end{array}%
\right) ,j=2,5,8,...,n-1  \label{17}
\end{equation}

\begin{equation}
\rho _{j}=\left( 
\begin{array}{c}
0 \\ 
0 \\ 
h_{1}U_{\frac{j-3}{3}}(\frac{\lambda _{1}}{2}) \\ 
0 \\ 
0 \\ 
h_{2}U_{\frac{j-3}{3}}(\frac{\lambda _{2}}{2}) \\ 
\vdots \\ 
0 \\ 
0 \\ 
h_{\frac{n}{3}}U_{\frac{j-3}{3}}(\frac{\lambda _{\frac{n}{3}}}{2})%
\end{array}%
\right) ,j=3,6,9,...,n-3,n  \label{19}
\end{equation}

From (\ref{17}) and\ (\ref{19}), the inverse of matrix $P$ can be written as
below;

\begin{eqnarray}
P^{-1} &=&\left[ 
\begin{array}{ccccccc}
h_{1}U_{0}(\frac{\lambda _{1}}{2}) & 0 & 0 & h_{1}U_{1}(\frac{\lambda _{1}}{2%
}) & 0 & 0 & \cdots \\ 
0 & h_{1}U_{0}(\frac{\lambda _{1}}{2}) & 0 & 0 & h_{1}U_{1}(\frac{\lambda
_{1}}{2}) & 0 & \cdots \\ 
0 & 0 & h_{1}U_{0}(\frac{\lambda _{1}}{2}) & 0 & 0 & h_{1}U_{1}(\frac{%
\lambda _{1}}{2}) & \cdots \\ 
h_{2}U_{0}(\frac{\lambda _{2}}{2}) & 0 & 0 & h_{2}U_{1}(\frac{\lambda _{2}}{2%
}) & 0 & 0 & \cdots \\ 
0 & h_{2}U_{0}(\frac{\lambda _{2}}{2}) & 0 & 0 & h_{2}U_{1}(\frac{\lambda
_{2}}{2}) & 0 & \cdots \\ 
0 & 0 & h_{2}U_{0}(\frac{\lambda _{2}}{2}) & 0 & 0 & h_{2}U_{1}(\frac{%
\lambda _{2}}{2}) & \cdots \\ 
\vdots & \vdots & \vdots & \vdots & \vdots & \vdots & \ddots \\ 
h_{\frac{n}{3}}U_{0}(\frac{\lambda _{\frac{n}{3}}}{2}) & 0 & 0 & h_{\frac{n}{%
3}}U_{1}(\frac{\lambda _{\frac{n}{3}}}{2}) & 0 & 0 & \cdots \\ 
0 & h_{\frac{n}{3}}U_{0}(\frac{\lambda _{\frac{n}{3}}}{2}) & 0 & 0 & h_{%
\frac{n}{3}}U_{1}(\frac{\lambda _{\frac{n}{3}}}{2}) & 0 & \cdots \\ 
0 & 0 & h_{\frac{n}{3}}U_{0}(\frac{\lambda _{\frac{n}{3}}}{2}) & 0 & 0 & h_{%
\frac{n}{3}}U_{1}(\frac{\lambda _{\frac{n}{3}}}{2}) & \cdots%
\end{array}%
\right.  \notag \\
&&\left. 
\begin{array}{cccc}
\cdots & h_{1}U_{\frac{n-3}{3}}(\frac{\lambda _{1}}{2}) & 0 & 0 \\ 
\cdots & 0 & h_{1}U_{\frac{n-3}{3}}(\frac{\lambda _{1}}{2}) & 0 \\ 
\cdots & 0 & 0 & h_{1}U_{\frac{n-3}{3}}(\frac{\lambda _{1}}{2}) \\ 
\cdots & h_{2}U_{\frac{n-3}{3}}(\frac{\lambda _{2}}{2}) & 0 & 0 \\ 
\cdots & 0 & h_{1}U_{\frac{n-3}{3}}(\frac{\lambda _{2}}{2}) & 0 \\ 
\cdots & 0 & 0 & h_{1}U_{\frac{n-3}{3}}(\frac{\lambda _{2}}{2}) \\ 
\ddots & \vdots & \vdots & \vdots \\ 
\cdots & h_{\frac{n}{3}}U_{\frac{n-3}{3}}(\frac{\lambda _{\frac{n}{3}}}{2})
& 0 & 0 \\ 
\cdots & 0 & h_{\frac{n}{3}}U_{\frac{n-3}{3}}(\frac{\lambda _{\frac{n}{3}}}{2%
}) & 0 \\ 
\cdots & 0 & 0 & h_{\frac{n}{3}}U_{\frac{n-3}{3}}(\frac{\lambda _{\frac{n}{3}%
}}{2})%
\end{array}%
\right] \text{.}  \label{+++}
\end{eqnarray}

Here, we obtain the $h_{k}$ related to eigenvalues of $H$ as the similar
meaning in [2-4].

\begin{equation}
h_{k}=\frac{3\left( -1\right) ^{k}\left( \lambda _{k}^{2}-4\right) }{2n+6}%
,n=3k\text{ }(k=2,4,6,...,n/3)  \label{20}
\end{equation}

\begin{equation}
h_{k}=\left\{ 
\begin{array}{c}
\frac{3(-1)^{k+1}\lambda _{\frac{n+6k+3}{6}}^{2}}{2n+6}\text{ },\text{ }if%
\text{ }1\leq k\leq \frac{n-3}{6} \\ 
\frac{6(-1)^{k+1}}{n+3}\text{ },\text{ }if\text{ }k=\frac{n+3}{6} \\ 
\frac{3(-1)^{k+1}\lambda _{\frac{n-2k+3}{2}}^{2}}{2n+6}\text{ },\text{ }if%
\text{ }\frac{n+9}{6}\leq k\leq \frac{n}{3}%
\end{array}%
\right. ,n=3k\ (k=3,5,7,...,n/3)  \label{21}
\end{equation}

By using (\ref{+}), (\ref{++}),(\ref{+++}) and (\ref{1}) we obtain general
expression as follow;

for $n=3k,$ $k=1,3,5,...,n/3,$ $k\in 
\mathbb{N}
$;

\begin{eqnarray}
\left\{ H^{m}\right\} _{i,j} &=&\dsum\limits_{k=1}^{\frac{n-3}{6}}\left(
\lambda _{2k-1}^{m}h_{2k-1}U_{\frac{i-\delta _{ij}}{3}}(\frac{\lambda _{2k-1}%
}{2})U_{\frac{j-\delta _{ij}}{3}}(\frac{\lambda _{2k-1}}{2})+\lambda
_{2k}^{m}h_{2k}U_{\frac{n-\sigma _{ij}}{3}}(\frac{\lambda _{2k}}{2})U_{\frac{%
j-\delta _{ij}}{3}}(\frac{\lambda _{2k}}{2})\right)  \notag \\
&&+\lambda _{\frac{n}{3}}^{m}h_{\frac{n}{3}}U_{\frac{i-\delta _{ij}}{3}}(%
\frac{\lambda _{\frac{n}{3}}}{2})U_{\frac{j-\delta _{ij}}{3}}(\frac{\lambda
_{\frac{n}{3}}}{2})  \label{28}
\end{eqnarray}%
where $h_{k}$ is defined as in (\ref{21}),

for $n=3k,k=2,4,6,...,n/3,k\in 
\mathbb{N}
;$

\begin{equation}
\left\{ H^{m}\right\} _{i,j}=\dsum\limits_{k=1}^{\frac{n}{6}}\left( \lambda
_{2k-1}^{m}h_{2k-1}U_{\frac{i-\delta _{ij}}{3}}(\frac{\lambda _{2k-1}}{2})U_{%
\frac{j-\delta _{ij}}{3}}(\frac{\lambda _{2k-1}}{2})+\lambda
_{2k}^{m}h_{2k}U_{\frac{n-\sigma _{ij}}{3}}(\frac{\lambda _{2k}}{2})U_{\frac{%
j-\delta _{ij}}{3}}(\frac{\lambda _{2k}}{2})\right)  \label{29}
\end{equation}%
where $h_{k}$ is defined as in (\ref{20}). $\delta _{ij}$ and $\sigma _{ij}$
are expressed as follows each formula;

\begin{equation*}
\delta _{ij}=\left\{ 
\begin{array}{c}
1\text{ \ \ \ \ },\text{ \ \ \ \ }i+j\equiv 2\func{mod}(3) \\ 
2\text{ \ \ \ \ },\text{ \ \ \ \ }i+j\equiv 1\func{mod}(3) \\ 
3\text{ \ \ \ \ },\text{ \ \ \ \ }i+j\equiv 0\func{mod}(3)%
\end{array}%
\right.
\end{equation*}

and

\begin{equation*}
\sigma _{ij}=\left\{ 
\begin{array}{c}
i+2\text{ \ \ \ \ },\text{ \ \ \ \ }i+j\equiv 2\func{mod}(3) \\ 
i+1\text{ \ \ \ \ },\text{ \ \ \ \ }i+j\equiv 1\func{mod}(3) \\ 
i\text{ \ \ \ \ \ \ \ \ \ },\text{ \ \ \ \ }i+j\equiv 0\func{mod}(3)%
\end{array}%
\right. \text{.}
\end{equation*}

\section{Numerical Considerations}

Let's give two numerical examples for the matrix (\ref{Hn}).

For$\ n=6,J=diag(\lambda _{1},\lambda _{1},\lambda _{1},\lambda _{2},\lambda
_{2},\lambda _{2})$ and from the general formula (\ref{29}), we obtain the
each elements of $\left\{ H^{m}\right\} _{ij}$ as

\begin{equation*}
h_{11}=h_{22}=h_{33}=\frac{1}{6}(\lambda _{1}^{m}(-1)^{1}(\lambda
_{1}^{2}-4)+\lambda _{2}^{m+1}(-1)^{2}(\lambda _{2}^{2}-4))
\end{equation*}%
\begin{equation*}
h_{14}=h_{25}=h_{36}=\frac{1}{6}(\lambda _{1}^{m+1}(-1)^{1}(\lambda
_{1}^{2}-4)+\lambda _{2}^{m+2}(-1)^{2}(\lambda _{2}^{2}-4))
\end{equation*}%
\begin{equation*}
h_{44}=h_{55}=h_{66}=\frac{1}{6}(\lambda _{1}^{m+2}(-1)^{1}(\lambda
_{1}^{2}-4)+\lambda _{2}^{m+1}(-1)^{2}(\lambda _{2}^{2}-4))
\end{equation*}%
\begin{equation*}
h_{41}=h_{52}=h_{63}=\frac{1}{6}(\lambda _{1}^{m+1}(-1)^{1}(\lambda
_{1}^{2}-4)+\lambda _{2}^{m}(-1)^{2}(\lambda _{2}^{2}-4))\text{.}
\end{equation*}

where $\lambda _{k}$ is as (\ref{11}) and also other components are zero.

For $n=9,$ let $J=diag(\lambda _{1},\lambda _{1},\lambda _{1},\lambda
_{2},\lambda _{2},\lambda _{2},\lambda _{3},\lambda _{3},\lambda _{3})$ and
from the general formula (\ref{28}), we obtain the each elements of $\left\{
H^{m}\right\} _{ij}$ as

\begin{equation*}
a_{11}=a_{22}=a_{33}=\frac{1}{8}[(\lambda _{1}^{m}+\lambda
_{3}^{m})(-1)^{2}\lambda _{3}^{2}+4\lambda _{2}^{m}(-1)^{2}(\lambda
_{2}^{2}-1)]
\end{equation*}%
\begin{equation*}
a_{14}=a_{25}=a_{36}=\frac{1}{8}[(\lambda _{1}^{m+1}+\lambda
_{3}^{m+1})(-1)^{2}\lambda _{3}^{2}+4\lambda _{2}^{m}(-1)^{2}(\lambda
_{2}^{3}-\lambda _{2})]
\end{equation*}%
\begin{equation*}
a_{17}=a_{28}=a_{39}=\frac{1}{8}[(\lambda _{1}^{m}(\lambda
_{1}^{2}-1)+\lambda _{3}^{m}(\lambda _{3}^{2}-1))(-1)^{2}\lambda
_{3}^{2}+4\lambda _{2}^{m}(-1)^{2}(\lambda _{2}^{2}-1)^{2}]
\end{equation*}%
\begin{equation*}
a_{41}=a_{52}=a_{63}=\frac{1}{8}[(\lambda _{1}^{m+1}+\lambda
_{3}^{m+1})(-1)^{2}\lambda _{3}^{2}+4\lambda _{2}^{m+1}(-1)^{2}]
\end{equation*}%
\begin{equation*}
a_{44}=a_{55}=a_{66}=\frac{1}{8}[(\lambda _{1}^{m+2}+\lambda
_{3}^{m+2})(-1)^{2}\lambda _{3}^{2}+4\lambda _{2}^{m+2}(-1)^{2}]
\end{equation*}%
\begin{equation*}
a_{47}=a_{58}=a_{69}=\frac{1}{8}[(\lambda _{1}^{m}(\lambda _{1}^{3}-\lambda
_{1})+\lambda _{3}^{m}(\lambda _{3}^{3}-\lambda _{3}))(-1)^{2}\lambda
_{3}^{2}+4\lambda _{2}^{m}(\lambda _{2}^{3}-\lambda _{2})]
\end{equation*}%
\begin{equation*}
a_{71}=a_{82}=a_{93}=\frac{1}{8}[(\lambda _{1}^{m}(\lambda
_{1}^{2}-1)+\lambda _{3}^{m}(\lambda _{3}^{2}-1))(-1)^{2}\lambda
_{3}^{2}+4\lambda _{2}^{m}(-1)^{2}]
\end{equation*}%
\begin{equation*}
a_{74}=a_{85}=a_{96}=\frac{1}{8}[(\lambda _{1}^{m}(\lambda _{1}^{3}-\lambda
_{1})+\lambda _{3}^{m}(\lambda _{3}^{3}-\lambda _{3}))(-1)^{2}\lambda
_{3}^{2}+4\lambda _{2}^{m+1}(-1)^{2}]
\end{equation*}%
\begin{equation*}
a_{77}=a_{88}=a_{99}=\frac{1}{8}[(\lambda _{1}^{m}(\lambda
_{1}^{2}-1)^{2}+\lambda _{3}^{m}(\lambda _{3}^{2}-1)^{2})(-1)^{2}\lambda
_{3}^{2}+4\lambda _{2}^{m}(-1)^{2}(\lambda _{2}^{2}-1)]\text{.}
\end{equation*}

where $\lambda _{k}$ is as (\ref{11}) and also other components are zero.

\section{Conclusion}

We investigate a symmetric Heptadiagonal Matrix in terms of positive integer
powers and relationship with second kind Chebyshev Polynomials. In this
context, we obtain the eigenvalue formula of (\ref{Hn}) associate with
second kind Chebyshev Polynomial roots. Additionally, we obtain the general
formula of positive integer powers for (\ref{Hn}). In the last section of
paper, some illustrative examples are given.

\end{document}